\documentclass[conference]{IEEEtran}
\IEEEoverridecommandlockouts

\usepackage{cite}
\usepackage{amsmath,amssymb,amsfonts}

\usepackage{algorithm}
\usepackage{graphicx}
\usepackage{textcomp}
\usepackage{xcolor}
\usepackage{subcaption}

\usepackage[bookmarks=false]{hyperref}

\usepackage{algorithm}
\usepackage{lipsum}

\makeatletter
\def\footnoterule{\kern-3\p@
  \hrule \@width 2in \kern 2.6\p@} 
\makeatother

\begin{document}

\title{Location Planning of Fast Charging Station considering its Impact on the Power Grid Assets}



\author{\IEEEauthorblockN{Daijiafan Mao\IEEEauthorrefmark{1}, Jun Tan\IEEEauthorrefmark{2}, Guangyi Liu\IEEEauthorrefmark{2}, and Jiankang Wang\IEEEauthorrefmark{1}} \IEEEauthorblockA{\IEEEauthorrefmark{1}Department of Electrical and Computer Engineering, The Ohio State University, Columbus, OH\\ Email: mao.156@osu.edu, wang.6536@osu.edu} \IEEEauthorblockA{\IEEEauthorrefmark{2}Global Energy Interconnection Research Institute North America, San Jose, CA\\ Email: jun.tan@geirina.net, guangyi.liu@geirina.net}}

\maketitle

\begin{abstract}
Under the ambition of boosting Plug-in Electric Vehicle (PEV) charging speed to a level comparable to the traditional refueling, Fast Charging Station (FCS) has been integrated into power distribution system. The location planning of FCS must allow for satisfactory charging service for PEV users as well as mitigate the detrimental effects on power grid caused by uncertainty and impulsiveness of charging demand. This paper proposed a location planning model for FCS, taking into account its impacts on the critical power grid assets. The multi-objective planning model simultaneously considered the role of FCS in the electricity and transportation sectors. This planning model is solved by the cross-entropy (CE) method. The validity and effectiveness of the CE approach have been demonstrated on a synthetic coupled network.
\end{abstract}


\section{Introduction} \label{intro}
\let\thefootnote\relax\footnote{This work was sponsored by the Ford Motor Company.}
In the recent decade, the electric power distribution system has been increasingly penetrated with Plug-in Electric Vehicles (PEV) \cite{IEAEVOutlook2017}. As an unconventional type of electric load, PEV consumes a large amount of power in a fast and discrete manner, which far exceeds the peak power demand for an average household in the U.S \cite{mao2017evaluating}. In continuing towards the goal of resembling traditional refueling stations, multiple charge ports per station are built, bringing total power capacity in a Fast Charging Station (FCS) to Megawatt level \cite{meyer2018integrating}. Prior works have shown that the impulsive characteristics of PEV charging load, especially when aggregated at the FCS, will result in negative impacts on the power grid, including disruptively varying voltage profiles along the feeder and lifetime depreciation of critical grid assets \cite{mu2014spatial}.

From electricity sector's point of view, the goal of FCS planning mainly centers on optimally placing FCS in the electric power network to achieve minimal investment, operation \& maintenance, electricity costs, etc. \cite{liu2013optimal,zhang2016integrated}. Furthermore, the capacity sizing and economic evaluation of on-site distributed generation and/or energy storage for the charging station have been considered in FCS planning paradigm as well \cite{yuan2017co,mao2015economic,negarestani2016optimal}. Yet it is important to note that there are additional parties/sectors that are interested in the planning of FCS locations. While the existing literature has considered electrical or traffic impacts of PEV charging stations, the problem facing municipal planners is multi-disciplinary in nature. 

The traffic and electrical considerations interact with each other in the framework of FCS planning through the following ways: (i) PEV users must have FCS that allows them to travel from origin to destination without exceeding their electric range, (ii) users must have desirable FCS available in order to fully facilitate the adoption of PEVs with little change to daily driving habits. The objective of FCS location planning from electrical and traffic point of views are conflicting in general. For example, though it may be beneficial to the power grid planner to allocate an FCS away from existing load center or close to the head of the feeder for grid loss minimization, this particular location is likely undesirable to a large amount of PEV drivers and is thus, sub-optimal for them. Therefore, by coupling power distribution and transportation network, this paper is able to take the PEV's impacts on both the power grid assets and traffic conditions into consideration.

The aforementioned considerations have indicated that the allocation of PEV fast charging stations is a multi-faceted, interdisciplinary problem, as FCS must allow for complete and efficient travel service as well as minimal adverse impacts on power distribution grid. Therefore, an integrated planning model need to be formulated by combining the view of FCS as an impulsive electric load and an extension of human behavior. This paper proposed a location planning model of FCS considering its impacts on the critical power grid assets, and thus reinforced the role that FCS plays in interacting with coupled electric and transportation network.
\section{FCS Location Planning Model} \label{proposed}
\subsection{Integrated Algorithm for Evaluating PEV's Impact on the Power Grid Assets}\label{GADM}
 With increasing PEV adoptions and improving FCS paradigms, it is critical for electric utilities to accurately quantify the impact of PEV charging on grid assets and plan for equipment replacement and infrastructure expansion accordingly, in order to ensure service reliability. Previous works on such impact assessment do not naturally fulfill utilities' needs of quantifying the long-term cost of critical assets. One of the most important reasons is that the time-series power flow (TSPF) analysis are taken in the form of annual average in the grid asset assessment \cite{rural2016guide}, which makes PEVs’ impulsive charging characteristics invisible. In other words, the load spikes caused by PEV charging can be easily averaged off in the assessment and shown harmless, while they could greatly reduce the lifetime of the grid assets in reality.

Therefore, an integrated algorithm for evaluating PEV's impact on the power grid asset has been proposed in the authors' preliminary work \cite{mao2019integrated}. It provides a convenient assessment through an integrated interface and is capable of capturing the inter-temporal response of grid assets. The advantageous features of the proposed algorithm are realized under a unified mathematical framework, where grid assets' depreciation model (GADM) are established and their re-casted Total Cost of Ownership (TCO) during any time span of interest are exemplified by a substation transformer as follow
\begin{align}\label{eq:mtcot}\begin{split}
TCO (t_1,t_2) = &L_T (t_1,t_2)\cdot C_o + CL\cdot A(t_1,t_2) \\
&+ LL \cdot B(s,t_1,t_2),
\end{split}
\end{align}
where the definition of parameters and modification of TCO method in utility practice have been detailed in Ref. \cite{mao2019integrated}. Noted that this paper does not address the transient response and voltage instability induced by PEV charging in the impact metrics \cite{huang2013quasi,mao2018impact}.

The aggregated charging profiles at the FCS are pre-processed through Monte-Carlo Simulation (MCS), which ensures accounting of random charging patterns over time and space. Therefore, the GADM combined with re-established TCO evaluation can accurately capture any overloading form. Distinct from simulation-based methods, this integrated algorithm is analytical, and thus greatly reduce the computation resources and data required for accurate assessment.
\subsection{Location Planning Formulation} \label{model}
A key consideration for FCS location planning, outside of power grid characteristics, is the attraction of volume of traffic flow within the region of interest. Traffic flow modeling is commonly employed to locate infrastructures such as gas stations, retail facilities, and for expansion planning. As FCS are, in the future, essentially gas stations dispensing electricity instead of gas, it would be helpful to adopt existing methodologies for siting. In particular, the utilization of flow-capturing model (FCM) will improve the planning feasibility and objectives for system operator as it assumes that the PEV charging demands are generated from vehicle flows and that the charging action occurs \textit{en route} from origin to destination \cite{hodgson1990flow,wu2017stochastic}. It should be noted that the proposed planning model in this paper does not consider posterior impacts of charging facilities. That is, the common driving rationales (e.g., shortest distances, fastest time, toll cost, etc.) determine the routes; the locations of charging facilities do not alter the routes but only capture the routes.

Previous works on FCS location planning mainly focused on the minimization of costs associated with FCS itself, i.e., installation, operation and maintenance, induced expansion cost, etc. From utility's planning point of view, it would be important to investigate and then mitigate the long-term cost of grid assets induced by FCS while maintaining a desirable public charging service. The existing studies on siting FCS rarely took on this perspective. Therefore, the multi-faceted location planning of FCS in this paper combines a view of FCS as an impulsive load to the power grid and an extension of user behavior in transportation sector. Both grid aspects (electrical responses and utility considerations) and transportation aspect are effectively integrated to create one multi-objective planning model as follow
\begin{align}
\underset{x,y}{\text{minimize}}& \quad TCO(x,y) - c\sum_{q\in Q} f_q y_q \label{evaluate}\\ 
\text{s.t.     } \begin{split} &\quad TCO(x,y) = L_T(x,y) \cdot C_o + CL\cdot A \\ & \qquad \qquad \qquad \quad ~~+ LL\cdot B(s(x,y)) \end{split} \\
\quad & \quad \sum_{k \in K}x_k = N_{FCS} \label{cons_fcs}\\
\quad & \quad x_k\in \{0,1\} \qquad \forall k\in K \\
\quad & \quad y_q\in \{0,1\} \qquad \forall q\in Q \label{yq},
\end{align}
where $TCO$ is the total cost of substation transformer during a time span (viz. equation (\ref{eq:mtcot})), $L_T(x)$ is the daily loss of life of transformer, $C_o$ is the capitol cost of transformer, $N_{FCS}$ is the total number of FCS installed within the distribution system, $Q$ is the set of PEV trip chain indices on the road network, $K$ is the set of all candidate FCS locations on the road network, $c$ is the coefficient that transforms the captured volume into the monetary cost, $f_q$ is the number of PEVs traveling on trip chain q, $x_k$ is the binary variable that equals 1 if a station is placed at location k and equals 0 otherwise, $y_q$ is the binary variable that equals 1 if PEVs on trip chain q are captured by a charging station and equals 0 otherwise.

As mentioned above, given a fixed origin-destination pair (O-D pair) for each individual PEV, the detailed paths, i.e., trip chain $q$, in the transportation network is pre-determined by a common rationale such as shortest distance route. Trip paths of all PEVs (set $Q$) are then evaluated in FCM to calculate the total captured traffic volume among all candidate FCS. The algorithm of determining whether or not a trip chain can be captured is adopted from Ref. \cite{wu2017stochastic}. The total captured PEV volume forms the temporal charging demands that impose on substation transformer, which cause more frequent loading spikes and adverse lifetime degradation of the assets, per the conclusion in author's previous work \cite{mao2019integrated}.

\begin{algorithm*}
\caption{CE Algorithm for Location Planning of FCS}
\label{ce_alg} 
Parameters:\\
$N_{FCS}$: number of FCS built in the area; $M=|\mathcal{X}|$: number of candidate locations of FCS; \\
$N=2000$: number of sample solutions in each iteration; $\rho=0.05$: rarity parameter; $\alpha=0.7$: smoothing parameter.\\\\
1: Initialize $\hat{\mathbf{v}}_0$ with Bernoulli density $\hat{\mathbf{v}}_{0,i} =N_{FCS}/M$ for $i=1,\ldots,M$. Set the iteration counter $t=0$.\\
2: Increment $t$ by 1. Generate $X_1,\ldots,X_N \sim_{i.i.d.} f(\cdot;\hat{\mathbf{v}}_{t-1})=\prod_{j=1}^n v_j^{x_j}(1-v_j)^{1-x_j}$. Let the performance function $S(\mathbf{X}_t)$ be the objective function in (\ref{evaluate}). Evaluate performance function values $S(X_j)$ for $j=1,\ldots,N$ and sort them from the smallest to largest such that $\{S(X_1)\leq \ldots \leq S(X_N)\}$. The solution vectors $\mathbf{X}_t$ are ordered accordingly.\\
3: After sorting, the best few solutions determined by rarity parameter $\rho$ are selected for updating the parameter vector $\mathbf{\hat{v}}$ as $\hat{\mathbf{v}}_{t+1}=\text{mean}(\mathbf{X}^{\text{elites}}_t)$, where $\mathbf{X}^{\text{elites}}_t=\{X(1),X(2),\ldots,X(\rho N)\}$, i.e., elite samples of $\mathbf{X}_t$.\\
4: To avoid the guided search getting trapped in local minimum and improve searching ability, the parameter vector $\mathbf{\hat{v}}$ is passed through the \textit{smooth updating}: $\mathbf{\hat{v}}_{t+1}=\alpha \hat{\mathbf{v}}_{t+1}+(1-\alpha)\hat{\mathbf{v}}_{t}$, where $\alpha \in [0,1]$ is the smoothing parameter which can be adjusted iteratively. \\
5: Repeat from Step 2 until the vector $\hat{\mathbf{v}}_t$ has converged to a degenerate, i.e., binary, vector.
\end{algorithm*}
\subsection{Solving Algorithm} \label{subsec_algorithm}
The planning model in \ref{model} is a generalization of Knapsack Problem (KP) \cite{kellerer2003knapsack}. There is strong theoretical evidence that the KP and hence its variants, i.e., combinatorial problems, belong to a class of NP-hard problem \cite{lam2014electric}. The addition of highly nonlinear constraints of grid asset aspect, i.e., the $TCO$ term, further complicated the search space of the problem. Hence, it can be seen that the proposed FCS location planning problem is at least as hard as KP. Therefore, the heuristic method, in particular, cross-entropy (CE) optimization method, has been employed to solve this planning problem.

The CE method uses the cross-entropy (or KL-divergence) as a measure of closeness between two sampling distributions. The basic idea is that locating an optimal or near-optimal solution through random search is a rare event. The CE method can gradually steer the sampling distribution of the random search so that the rare event is more likely to occur \cite{rubinstein2001combinatorial}. Therefore, the CE method unifies many existing population-based optimization heuristics. 

The algorithm involves an iterative procedure where each iteration can be broken down into two steps \cite{de2005tutorial}. Suppose the search space is the finite set $\mathcal{X}$ with cardinality $|\mathcal{X}|=n$. In the first step, the deterministic problem has been randomized by defining a family of probability density functions (pdf) to generate sample objects $\mathbf{X}$ in $\mathcal{X}$ with sample size $N$. Since the solution $\mathbf{x}$ in (\ref{evaluate}) is a binary vector, a simple choice for the sampling pdf could be the multivariate Bernoulli density. Then in the second step the updating rule of parameter vector $\hat{\mathbf{v}}$ at $t$-th iteration becomes
\begin{align} \label{update}
    \hat{\mathbf{v}}_{t,j}=\frac{\sum_{k=1}^N \mathbf{I}_{\{\hat{S}(X_k)\leq \hat{\gamma}_t\}}X_{k,j}}{\sum_{k=1}^N \mathbf{I}_{\{\hat{S}(X_k)\leq \hat{\gamma}_t\}}}, \quad j=1,\ldots,n
\end{align}
where $X_1,\ldots,X_N$ are i.i.d. copies of $\mathbf{X}\sim \{\hat{v}_{t-1,j}\}$, $X_{k,j}$ is the $j$-th component of the $k$-th random binary vector $X_k$ and $\mathbf{I}_{\{\cdot\}}$ is the indicator function. The set of vectors satisfying $S(X_k)\leq \hat{\gamma}_t$ is called the \textit{elite samples} generally determined by the \textit{rarity parameter} $\rho$ at iteration $t$. 

Note that the updating rule (\ref{update}) is the analytical solution of the cross-entropy minimization in (\ref{ce_min}), since the sampling pdf $\mathbf{v}$ is chosen such that it belongs to a Natural Exponential Family, i.e., Bernoulli density \cite{rubinstein2004combinatorial}.
\begin{align} \label{ce_min}
    \underset{\mathbf{v}}{\text{max}}&~ \frac{1}{N}\sum_{k=1}^N \mathbf{I}_{\{S(X_k) \leq \hat{\gamma}_t)\}}\ln f(X_k ; \mathbf{v}) 
\end{align}

The detailed CE procedure at each iteration for location planning of FCS is elaborated in Algorithm \ref{ce_alg}.

\subsection{Workflow} \label{subsec_workflow}
\begin{figure}[h!]
	\centering
	\includegraphics[scale=0.34]{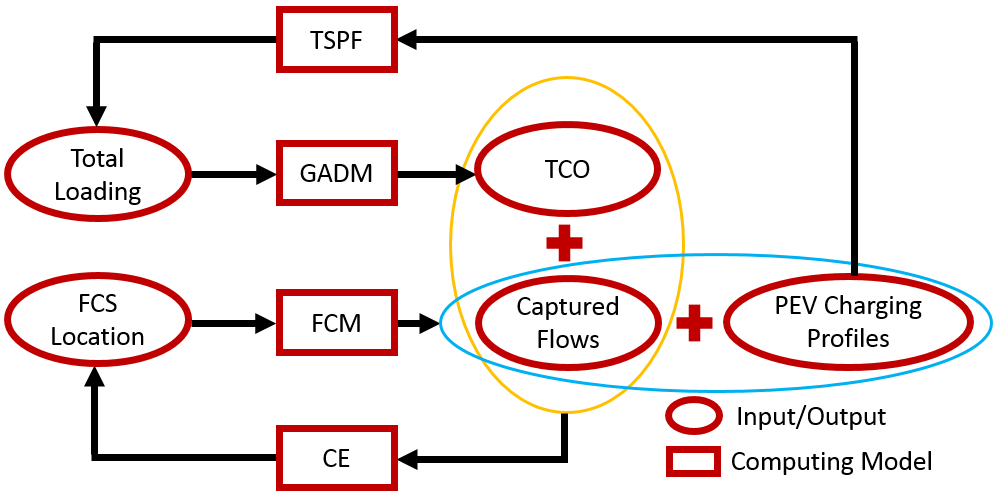}
	\caption{Workflow of Location Planning Model. TSPF: Time-Series Power Flow. GADM: Grid Asset Depreciation Model. FCM: Flow Capturing Model. CE: Cross-Entropy.}
	\label{workflow}
\end{figure}

The proposed location planning model and solving implementation has been outlined in Fig.\ref{workflow}. All computing models, i.e., TSPF, GADM, FCM and CE have been described in Section \ref{GADM} through Section \ref{subsec_algorithm}. The location of FCS will affect the daily captured PEV flows, which in turn alters the loading patterns at the corresponding electric node, and thus impact the substation transformer's TCO value. By optimizing the FCS locations, we can reduce PEVs' impact on grid assets while increase the FCS's charging service at the same time.
\section{Case Study} \label{case}
\subsection{Simulation Setup} \label{sim_setup}
\begin{figure*}[h!]
\begin{subfigure}{0.45\textwidth}
\includegraphics[width=0.85\linewidth, height=5.5cm]{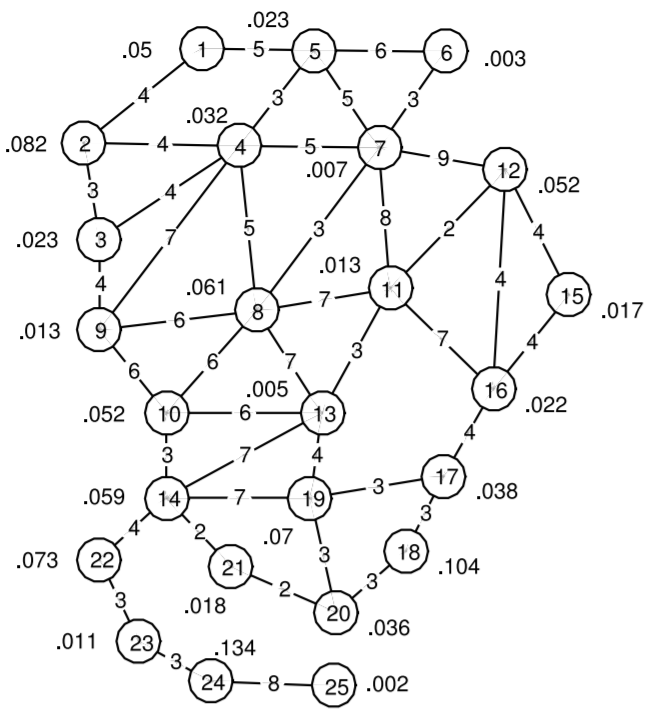} 
\caption{25-node Transportation Network}
\label{transportation}
\end{subfigure}
\begin{subfigure}{0.49\textwidth}
\includegraphics[width=1\linewidth, height=5.5cm]{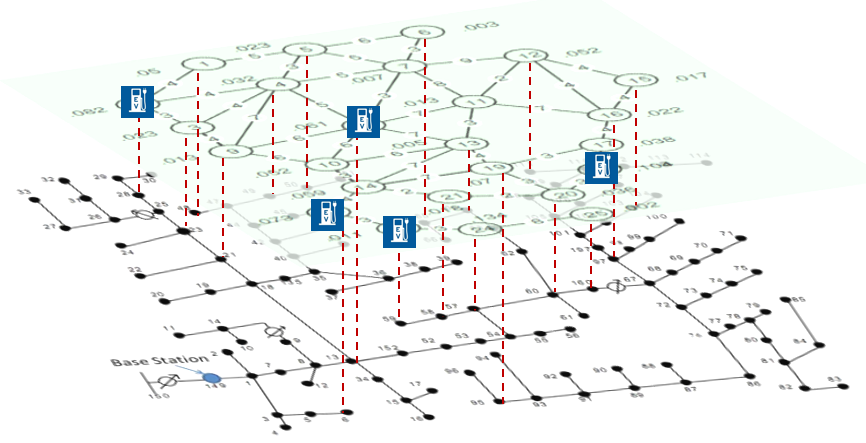} 
\caption{Illustration of Coupled Power and Transportation Network}
\label{coupled}
\end{subfigure}
\caption{Synthetic Test Network}
\label{network}
\end{figure*}
As the benchmark, the proposed location planning model with CE solving algorithm is implemented in MATLAB and tested on a synthetic network, i.e., 25-node transportation network coupled with IEEE 123-node test feeder, both of which are commonly analyzed in corresponding research communities \cite{kersting2001radial,hodgson1990flow}. All transportation nodes are randomly mapped to IEEE 123-node electric system, i.e., the transportation node set $\mathcal{N}^T$ is a subset of electric node set $\mathcal{N}^E$. The topology of transportation network and the electric-transportation coupled network are shown in Fig. \ref{transportation} and Fig. \ref{coupled}, respectively \cite{zhang2018pev}. The number on each transportation arc represents the normalized distance between the corresponding two nodes. We assume that the per-unit distance in Fig. \ref{transportation} is 5 mile, e.g., the distance of arc (1,2) is 4 units, which corresponds to 20 miles. 

\begin{figure}[h!]
	\centering
    \vspace{-5pt}
	\includegraphics[scale=0.46]{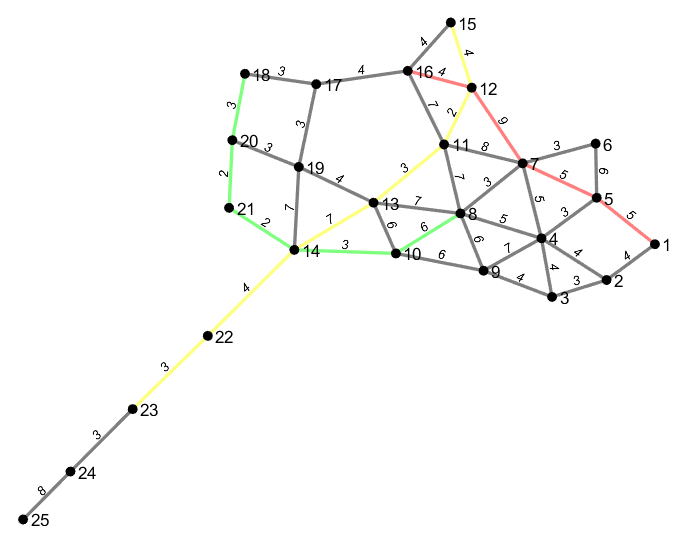}
	\caption{Illustration of Shortest-path Trip Chain for 3 PEVs}
	\label{path_illu}
\end{figure}

From an aggregation point of view, it is assumed that there is a total of 500 PEVs in this system, each of which has a random daily departure time, O-D pair and trip chain $q$ determined by shortest-path algorithm. Fig. \ref{path_illu} shows the MATLAB illustration of 3 trip chains in the 25-node transportation network, highlighted by red, green and yellow. The respective O-D pairs are (1,16), (8,18) and (15,23).

\subsection{Location Planning Results} \label{sim_results}
The performance function $S(\mathbf{X})$ in CE algorithm is equivalent to the objective function (\ref{evaluate}) of the planning model. Hence, the convergence of $S(\mathbf{X})$ simultaneously indicates the minimization of objectives. The constraint in (\ref{cons_fcs}), i.e., the resource capacity $N_{FCS}$ constraint, is enforced by adding a penalty term with coefficient $c_p=100$ in the objective function.

\begin{figure*}[h!]
\begin{subfigure}{0.49\textwidth}
\includegraphics[width=0.95\linewidth, height=6.3cm]{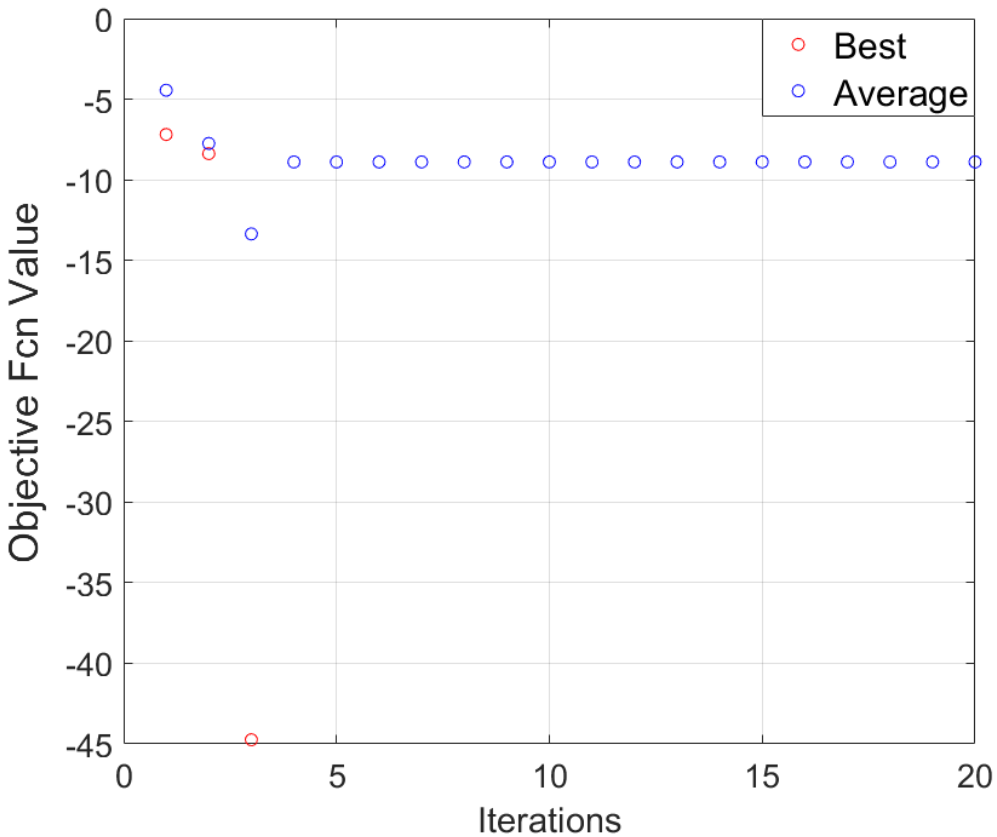} 
\caption{Convergence of Objective Function Value}
\label{conv_val_fcs5}
\end{subfigure}
\begin{subfigure}{0.49\textwidth}
\includegraphics[width=0.95\linewidth, height=6.3cm]{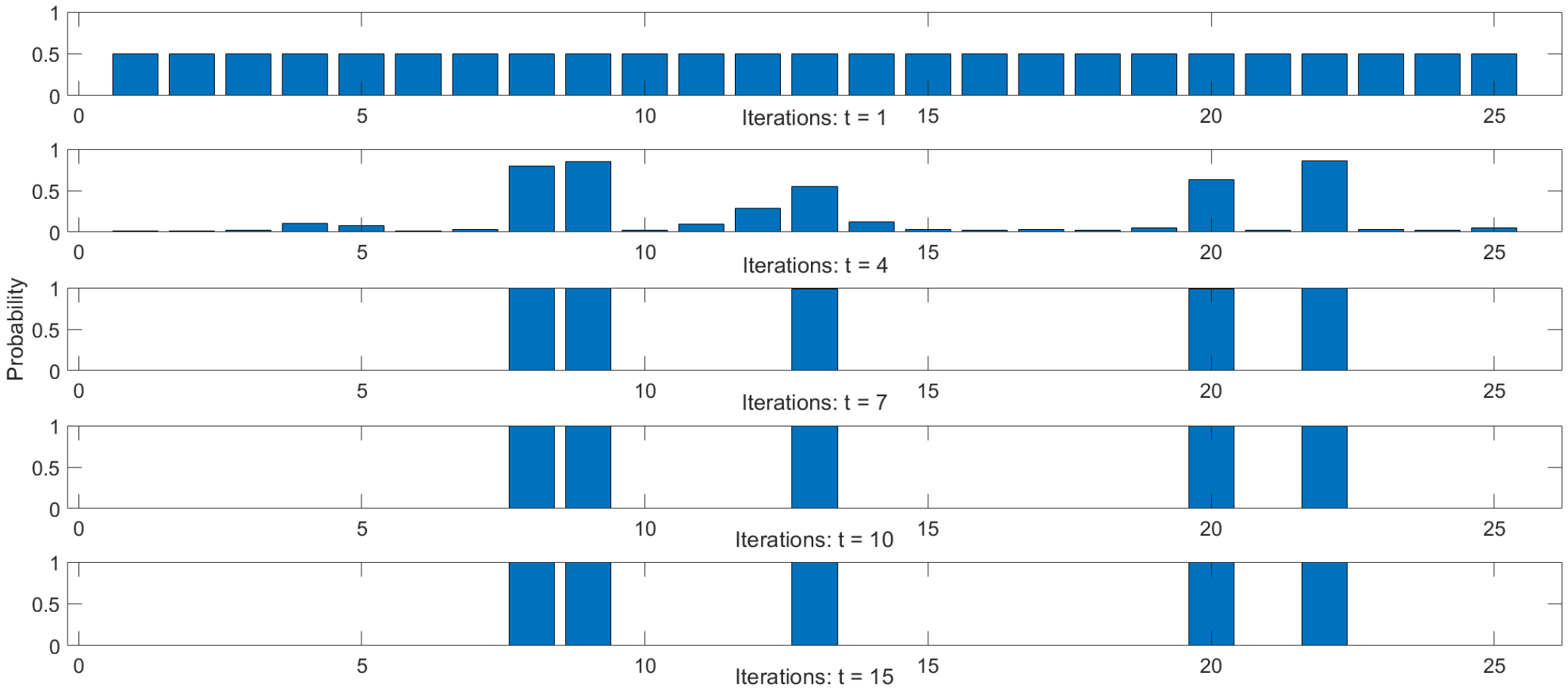} 
\caption{Evolution of Decision Variable Vector}
\label{conv_sol_fcs5}
\end{subfigure}
\caption{Location Planning Results for $N_{FCS}=5$}
\label{result_fcs5}
\end{figure*}

\begin{figure*}[h!]
\begin{subfigure}{0.49\textwidth}
\includegraphics[width=0.95\linewidth, height=6.3cm]{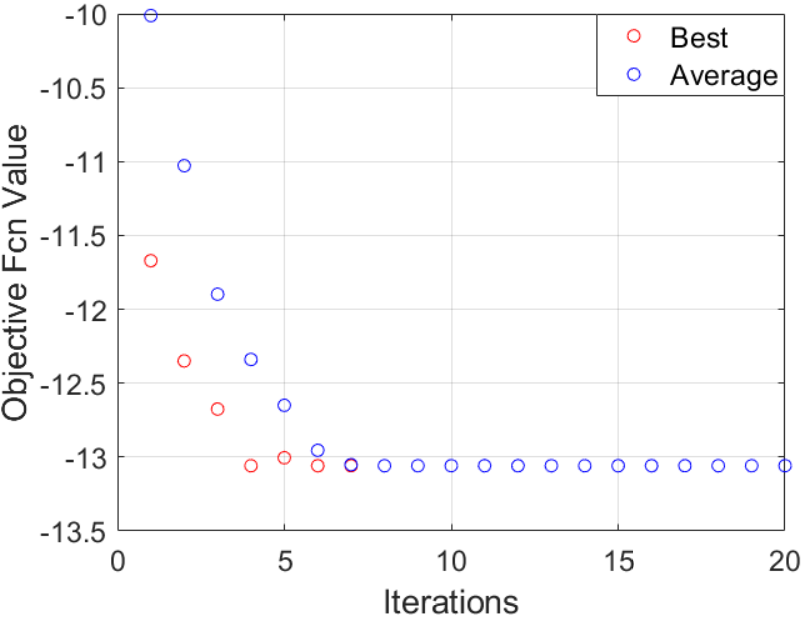} 
\caption{Convergence of Objective Function Value}
\label{conv_val_fcs10}
\end{subfigure}
\begin{subfigure}{0.49\textwidth}
\includegraphics[width=0.95\linewidth, height=6.3cm]{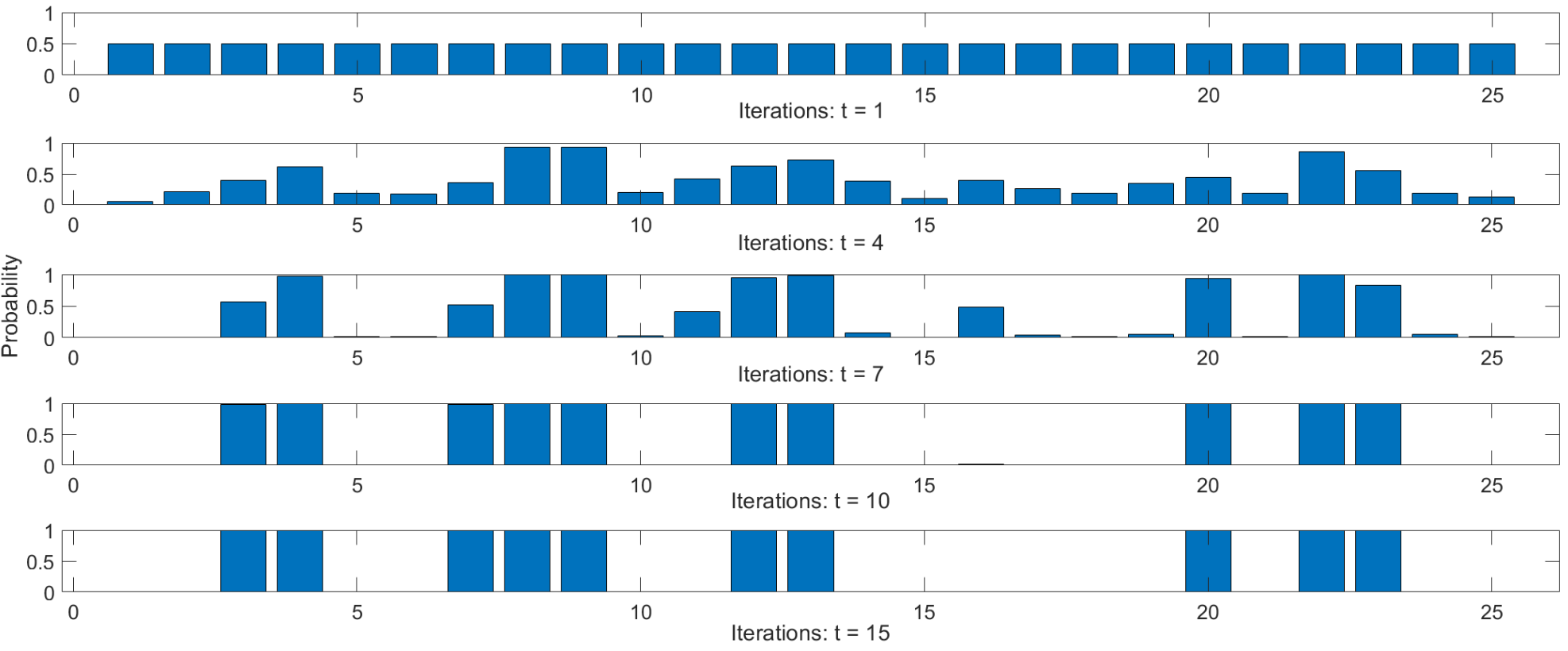} 
\caption{Evolution of Decision Variable Vector}
\label{conv_sol_fcs10}
\end{subfigure}
\caption{Location Planning Results for $N_{FCS}=10$}
\label{result_fcs10}
\vspace{-10pt}
\end{figure*}

The case study results solved by CE algorithm are presented in Fig. \ref{result_fcs5} and Fig. \ref{result_fcs10}, with the number of FCS built $N_{FCS}$ to be 5 and 10, respectively. As observed in Fig. \ref{conv_val_fcs5} and \ref{conv_val_fcs10}, the total cost, i.e., objective function value, converges in less than 20 CE iterations. Furthermore, Fig. \ref{conv_sol_fcs5} and \ref{conv_sol_fcs10} shows the evolution of probability vector $\hat{\mathbf{v}}_t$, which converges to a binary vector corresponding to the optimal solution. Therefore, the optimal locations to place FCS for $N_{FCS}=5$ and $N_{FCS}=10$ are $X = [8,9,13,20,22]$ and $X = [3,4,7,8,9,12,13,20,22,23]$, respectively.

In our future work, the in-depth algorithm verification and comparative study with other meta-heuristic algorithms (e.g., genetic algorithm) on the graph-computing platform will be conducted using a set of real-world power grid and transportation data.

\section{Conclusion} \label{conclusion}
This paper presented a location planning model for PEV fast charging stations (FCS), taking into account their impacts on the power grid assets. The multi-objective planning model integrally considered the role of FCS in both power and transportation sector, thus provide a more accurate assessing metric for utility planner. The proposed planning problem has been solved by the cross-entropy (CE) optimization method and verified on a synthetic power-transportation coupled network.

\bibliographystyle{ieeetran}
\bibliography{references}

\end{document}